\documentclass[12]{article}
\usepackage[centertags]{amsmath}
\usepackage{amsfonts}
\usepackage{amssymb}
\usepackage{amsthm}
\usepackage{newlfont}
\usepackage[ansinew]{inputenc}
\usepackage[english,francais]{babel}
\usepackage{amstext}

\theoremstyle{plain}
\newtheorem{Theorem}{Theorem}[section]
\newtheorem{proposition}{Proposition}[section]
\newtheorem{Corollary}{Corollary}[section]
\newtheorem{Definition}{Definition}[section]
\newtheorem{Lemma}{Lemma}[section]
\newtheorem{Remark}{Remark}[section]

\begin{document}

\title{Endomorphism Rings and Isogenies Classes for Drinfeld $A$-Modules of Rank 2 over Finite Fields }
\author{{\normalsize {MOHAMED AHMED Mohamed Saadbouh}} \\
{\small {Institut de Mathématiques de Luminy, CNRS-UPR 9016} }\\
{\small {Case 907, 163 Avenue de Luminy, 13288 Marseille, Cedex
9. France.}}}
\date{}
\maketitle

\selectlanguage{francais}
\begin{abstract}
Soit $\Phi $ un $\mathbf{F}_{q}[T]$-module de Drinfeld de rang
$2$, sur un corps fini $L$ $\mathbf{F}_{q}$, une extension de
degr\'{e} $n$ d'un corps fini On abordera plusieurs points
d'analogie avec les courbes elliptiques. Nous specifions les
conditons de maximalite et de non maximalite pour l'anneau
d'endomorphismes End$_{L}\Phi$
en tant que $\mathbf{F}_{q}[T]$-ordre dans l'anneau de division End$_{L}\Phi \otimes _{\mathbf{F}_{q}[T]}%
\mathbf{F}_{q}(T)$, on s'int\'{e}ressera ensuite aux polyn\^{o}me caract\'{e}%
ristique  et par son interm\'{e}diaire on calculera le nombre de
classes d'iog\'{e}nies.
\end{abstract}

\selectlanguage{english}

\begin{abstract}
Let $\Phi $ be a  Drinfeld  $\mathbf{F}_{q}[T]$-module of rank
$2$, over a finite field $L$, a finite extension of  $n$ degrees
of a finite field with $q$ elements $\mathbf{F}_{q}$. Let $m$ be
the extension degrees  of $ L$ over the field  $\mathbf{F}_{q}[T]/P$, $P$ is the $\mathbf{F}%
_{q}[T]$-characteristic of $L$, and $d$ the degree of the
polynomial $P$. We will discuss about a many analogies points
with elliptic curves. We start by the endomorphism ring of a
Drinfeld  $\mathbf{F}_{q}[T]$-module of rank 2, End$_{L}\Phi $,\
and we specify the maximality conditions and non maximality
conditions as a $\mathbf{F}_{q}[T]$-order in the ring of division End$_{L}\Phi \otimes _{\mathbf{F}_{q}[T]}%
\mathbf{F}_{q}(T)$, in the next point we will interest to the
characteristic polynomial of a Drinfeld module of rank 2 and used
it to calculate the number of isogeny classes for such module, at
last we will interested to the Characteristic of Euler-Poincare
$\chi _{\Phi }$ and we will calculated the cardinal of this
ideals.
\end{abstract}

\section{Introduction}

Let $E$ be a elliptic curves over a finite field
$\mathbf{F}_{q}$, we know, by [9], [12], [13] and [14], that the
endomorphism ring of $E$, End $_{\mathbf{F}_{q}}E$, is a  $
\mathbb{Z}-$order in a division algebra which is  : $ \mathbb{Q}$
and in this case  End$_{\mathbf{F}_{q}}E=\mathbb{Z}$, or is
complex quadratic field and in this case :
End$_{\mathbf{F}_{q}}E=\mathbb{Z}+c$ $O_{K}$ where $c$ is an
element of $ \mathbb{Z}$ and $O_{K}$ is the maximal
$\mathbb{Z}-$order in complex quadratic field maximal, or this
division algebra is a quaternion algebra $\mathbb{Q}$ in which
case End$_{F_{q}}E$ is a maximal order in this quaternion
algebra. We put  $E(\mathbf{F}_{q})$ the abelian group of
$\mathbf{F}_{q}$-rational  points of $E$. The cardinal of this
abelian group is equal to $N=q+1-c$, and by  Hasse-Weil $ \mid
c\mid \leq 2.\sqrt{q}$. A morphism of elliptic curve  over $
\mathbf{F}_{q}$, is a algebraic application $f:E_{1}\mapsto
E_{2}$, defined over $\mathbf{F}_{q}$, which respects the law
group. An isogeny is a no null morphism. And for one elliptic
curve $E$ the set of such morphism $f:E\mapsto E$ forms a ring  :
the  $\mathbf{F}_{q}$-endomorphism ring of $E$, this ring will be
noted End$_{\mathbf{F}_{q}}(E)$, of same the
$\bar{{\mathbf{F}_{q}}}$-endomorphism ring of $E$ will be noted
End$_{\overset{\_}{\mathbf{F}_{q}}}(E)$, in the case where
 End$_{\overset{\_}{\mathbf{F}_{q}}}(E)$ is non commutative the curve $E$ is called
 supersingular, otherwise it is called ordinary. And by [10], [13], [14]:

\begin{Theorem}
There are three possibilities for endomorphism ring for an\\
elliptic curve $E$ :

\begin{enumerate}
\item End$_{\mathbf{F}_{q}}(E)=\mathbb{Z}$

\item End$_{\mathbf{F}_{q}}(E)=\mathbb{Z}
+ c.O_{\max }$, where $c$ $\in \mathbb{Z}_{>0},$ $\ p$ not divide
$c$, and $O_{\max }$ is a maximal order in complex quadratic
field ( which is equal to the field of fractions of  End$
_{\mathbf{F}_{q}}(E)$( $c$ is called the conductor of
End$_{\mathbf{F}_{q}}(E)$)$;$

\item End$_{\mathbf{F}_{q}}(E)$ is a maximal order in the quaternion algebra $\mathbb{Q}_{\infty ,p}.$
\end{enumerate}
\end{Theorem}

The cardinal of  $E(\mathbf{F}_{q})$ and endomorphism Frobenius
$\varphi $ are given, in [9], [12] and [13], by :

\begin{Theorem}
 Let $E$ be an elliptic curve over  $\mathbf{F}_{q}$,  $
\varphi $ be the endomorphism Frobenius in
End$_{\mathbf{F}_{q}}(E)$ and $ p $ the  characteristic of
$\mathbf{F}_{q}:$

\begin{enumerate}
\item The endomorphism $\varphi $ satisfies an unique equation $
\varphi ^{2}-c\varphi +q=0$ in End$_{\mathbf{F}_{q}}(E)$, where $\
c\in \mathbb{Z}\subset $ End$_{\mathbf{F}_{q}}(E)$,

\item $\mid c\mid \leq 2\sqrt{q}$,

\item $\mid E(\mathbf{F}_{q})\mid =q+1-c$,

\item $p\mid c$ if and only if $E$ is supersingular.
\end{enumerate}
\end{Theorem}

The set of isogenies classes for an elliptic curve $E$ over a
finite field  $\mathbf{F}_{q}$, is given in [10], [13] and [14],
by :

\begin{Theorem}

The set of isogenies classes of an elliptic curve over a finite
field $\mathbf{F}_{q}$ is given by a natural bijection between
the set of integer  $c$ such that  $\mid c\mid \leq 2.\sqrt{q}$
and one of the following conditions is satisfied:

\begin{enumerate}
\item $(c,q)=1$;

\item $q$ is a square and $c=\pm 2\sqrt{q}$;

\item $q$ is a square, $p$ is not  congru to $1 (\mod 3)$, and$
c=\pm \sqrt{q}$;

\item $q$ is not a  square, $p=2$ or $3$, and $c=\pm \sqrt{p.q}$;

\item $q$ is not a square and $c=0$; or $q$ is a square, $p$ not
congru to  $1(\mod 4),$ and $c=0$.
\end{enumerate}
\end{Theorem}

The case 1 correspond to ordinary case and the other cases
correspond  to the supersingulary case.

Our goal here is to give an analog of the previous results in the
case of a Drinfeld Module of rank 2.We recall quickly what is it:
let $K$ a no empty global field of characteristic $p$ ( that
means a rational functions field of one indeterminate over a
finite field ) with a constant field the finite field
$\mathbf{F}_{q}$ with $p^{s}$ elements. We fix one place of  $K$,
noted $\infty $ \ and we call  $A$ \ the ring of regular elements
away from the place $\infty $. Let $L$ be a commutator field of
characteristic $p$, and let $\gamma :A\rightarrow L,$ be an
$A$-ring homomorphism,  the kernel of this  homomorphism is noted $P$ and $m$\ =$[L,$ $%
A/P]$ is the extension degrees of  $L$ over $A/P$.

We note $L\{\tau \}$ the Ore  polynomial ring, that means the
polynomial ring of  $\tau $, $\tau $ is  the Frobenius of
$\mathbf{F}_{q}$, with the usually  addition and the product is
given by the commutation rule : for every   $\lambda $ of
$L\mathit{,}$ $\tau \lambda =\lambda ^{q}\tau $. We said a
A-Drinfeld module $\Phi $ for a no trivial ring homomorphism,
from $A$ to \ $L\{\tau \}$ such that is different of $\gamma $.
this homomorphism $\Phi $, once defined, given a $A$-module
structure over the  $A$-field $L$, noted \ $L^{\Phi }$, where the
name of a Drinfeld $A$-module for a homomorphism $\Phi $. This
structure of $A$-module is depending on $\Phi $ and especially on
this rank.

Let $\chi $ be the characteristic of Euler-Poincaré ( it is a
ideal from $A$ ), so we can speak about the ideal  $\chi (L^{\Phi
})$, will be noted by $\chi _{\Phi }$, which is by definition a
divisor for $A$, corresponding for the elliptic curves to a number
of points of the variety over their  basic field.
We will work, in this paper, in the special case $K=$ $\mathbf{F}_{q}(T)$, $\ A=%
\mathbf{F}_{q}[T]$.
Let $P_{\Phi }(X):$ be the characteristic polynomial of the $A$%
-module $\Phi $, it is also a characteristic polynomial of
Frobenius $F$ of $L$.\ We can prove that this polynomial can be
given by : $P_{\Phi }(X)=$ $X^{2}-cX+\mu P^{m},$ such that $\mu
\in \mathbf{F}_{q}^{\ast } $, $c$ $\in A$ and $\ \deg c\leq
\frac{m.d}{2}$, the Hasse-Weil analogue in this case. We will
interested to the endomorphism ring of a Drinfeld  $A$-module of
rank 2 and we will prove:

\begin{proposition}
Let $\Delta $ $=c^{2}-4\mu P^{m},$ be the discriminant of $P_{F}$,
the  characteristic polynomial of the endomorphism Frobenius $F$
of the finite filed $L$, which is  $ P_{F}(X)=X^{2}-cX+\mu
P^{m}$, \ And let  $O_{K(F})$ the  maximal $A$-order of the
algebra $K(F)$.

\begin{enumerate}
\item For every $g\in A$ such that  $\Delta =g^{2}$ $.\omega $, there exists a Drinfeld $A$%
-module $\Phi $ over  $L$ of rank 2, such that : $O_{K(F})=A[\sqrt{%
\omega }]$ and :
\begin{equation*}
\text{End}_{L}\Phi =A+g.O_{K(F)}.
\end{equation*}

\item If there are not a polynomial $g$ of $A$ such that $g^{2}$ divide $%
\Delta $, then there exists a Drinfeld $A$-module  $\Phi $ over
$L$ of rank 2, such that End$_{L}\Phi =O_{K(F)}$.
\end{enumerate}
\end{proposition}
At next We will proved  that the number of isogenies classes of a
Drinfeld  $A$-module of rank 2 is :

\begin{proposition}
 Let $\Phi $ be a Drinfeld $A$-module of  rank $2$ over a finite filed $L=F_{q^{n}}$
 and let $P$ the $A$-characteristic of $L$. We put
 $m=[L:A/P]$ and $d=$deg $P$ :

\begin{enumerate}
\item  $m$ is odd and $d$ odd:%
\begin{equation*}
\#\{P_{\Phi }\text{, }\Phi :\text{ordinary}(1)\}=(q-1)(q^{[\frac{m}{2}%
d]+1}-q^{[\frac{m-2}{2}d]+1}+1).
\end{equation*}

\item $m$ is even and $d$ is odd :%
\begin{equation*}
\#\{P_{\Phi
}\}=(q-1)[\frac{q-1}{2}q^{\frac{m}{2}d}-q^{\frac{m-2}{2}d+1}+q].
\end{equation*}

\item $m$ is even and $d$ is even:%
\begin{equation*}
\#\{P_{\Phi
}\}=(q-1)[\frac{q-1}{2}q^{\frac{m.}{2}d}-q^{\frac{m-2}{2}d}+1].
\end{equation*}
\end{enumerate}
\end{proposition}

\section{Drinfeld Modules }

Let $E$ be an extension of $\mathbf{F}_{q}$, and let $\tau $ be the Frobenius of $%
\mathbf{F}_{q}$. We put $E\{\tau \}$ the polynomial ring in $\tau
$ with the usual  addition and the multiplication defined by:$\ $

\begin{equation*}
\ \forall e\in E,\text{ }\tau e=e^{q}\tau \text{.}
\end{equation*}
\begin{Definition}

Let  $R$ be  the $E$-linearly polynomials set with the
coefficient in
$E$, that means that these elements are on the following form :%
\begin{equation*}
Q(x)=\sum_{K>0}l_{k}x^{q^{K}},
\end{equation*}%
where $l_{k}\in E$ for every $k>0$, only a finite number of
$l_{k}$ is not null. The ring  $R$ is a ring by addition and the
polynomials composition.
\end{Definition}

\begin{Lemma}
$E\{\tau \}$ and $R$ are isomorphic rings.
\end{Lemma}

If we put  \ $A=\mathbf{F}_{q}[T],$ $\ f(\tau
)=\sum_{i=0}^{v}a_{i}\tau ^{i}\in E\mathit{\{\tau \}}$ \textit{and
\ }$Df\mathit{:=a}_{0}=f^{\prime }(\tau ).$

It is clear that the application :
\begin{equation*}
E\mathit{\{\tau \}\mapsto }E
\end{equation*}

\begin{equation*}
f\mapsto Df,
\end{equation*}%
is a  morphism of $\mathbf{F}_{q}$-algebra.
\begin{Definition}
\bigskip An $A$-field $E$  is a field $\mathit{E}$
equipped with a fix morphism \\ $\gamma :A \longrightarrow E$.
The prime ideal $P=$ Ker $\gamma $ is called the
$A$-characteristic of $E$.\\ We say $E$ has generic
characteristic if and only if $P=(0)$; otherwise
 (i.e $P\neq (0) $ ) we said  $E$ has a finite characteristic.
\end{Definition}

Then we have the following fundamental Definition :

\begin{Definition}
Let $E$ be an $A$-field and let $\Phi :A\mapsto E\mathit{\{\tau
\}}$ be a homomorphism of algebra. Then $\Phi $ is a  $A$-Drinfeld
module over $E$ if and only if:

\begin{enumerate}
\item $D\circ \Phi =\gamma ;$

\item for some $a$ $\in A,\Phi _{a}\neq \gamma (a)\tau ^{0}.$
\end{enumerate}
\end{Definition}
\begin{Remark}
\begin{enumerate}
\item the  normalization above is analogous  to the normalization used in complex
multiplication of elliptic curves. The last condition is
obviously a non-triviality condition.
\item By $\Phi $, every  extension $E_{0}$ of $E$ became an  $A$-module
by :
\begin{eqnarray*}
E_{0}\times A &\rightarrow &E_{0}, \\
\ (k,a) &\rightarrow &k.a:=\Phi _{a}(k)\ .
\end{eqnarray*}%
We will note this  $A$-module by $E_{0}^{\Phi }.$
\end{enumerate}
\end{Remark}

Let $\overline{E}$ a fix algebraic closure of $E$ and $\Phi $ a
Drinfeld module over $E$ and $I$ an ideal of $A$. As $A$ is a
Dedekind domain, one know that $I$ may be generated by ( at most ) two elements $%
\{i_{1}$ , $i_{2}\}\subset I.$

Since $E\mathit{\{\tau \}}$ has a right division algorithm, there
exists a right greatest common divisor in $E\mathit{\{\tau \}}$.
It is the monic generator of the left ideal of $E\mathit{\{\tau
\}}$ generated by : $\Phi _{i_{1}}$ et $\Phi _{i_{2}}.$

\begin{Definition}
We set  $\Phi _{I}$ to be the monic generator of the left ideal
of $E\mathit{\{\tau \}}$ generated by   $\Phi _{i_{1}}$ and
$\Phi_{i_{2}}.$
\end{Definition}

\begin{Definition}
Let $E_{0}$ be an extension of  $E$ and $I$ an ideal of $A$. We
define by  $: $ $\ \Phi \lbrack I](E)$ the finite subgroup of $\Phi \lbrack \overline{%
E_{0}}]$ given by the roots of $\Phi _{I}$ in $\overline{E}$.
\end{Definition} If  $a\in A$, then we set $\Phi \lbrack a]:=\Phi
\lbrack (a)].$ We can see it as :
 $\Phi \lbrack a]=\{$ set of roots of $\Phi (a)$ in $\overline{K}%
\},$ and $\Phi _{I}=\cap _{a\in I}\Phi \lbrack a].$ Then :

\begin{equation*}
\Phi _{a}(\overline{E}):=\Phi \lbrack a](\overline{E})=\{x\in \overline{E}%
,\Phi _{a}(x)=0\}
\end{equation*}

For every ideal  $Q\subset A,$
\begin{equation*}
\Phi _{Q}(\overline{E}):=\Phi _{Q}(\overline{E})=\cap _{a\in Q}\Phi _{a}(%
\overline{E}).
\end{equation*}

\begin{Remark}
The groups \ : $\Phi \lbrack I](E)$ and $\Phi \lbrack
I](\overline{E})$ are clearly stable under  $\{\Phi _{a}\}_{a\in
A}$.\ \ \
\end{Remark}

\begin{Definition}
Let $\Phi $ be a Drinfeld $A$-module over an  $A$-field $E$. We
say that  $ \Phi $ is supersingular, if and only if, the
$A$-module constituted by a $P$-division points $\Phi
_{P}(\overline{E})$ is trivial.
\end{Definition}


\subsection{The Height and Rank of a Drinfeld Module $\Phi $}

Let $\Phi $ be a Drinfeld $A$-module over the $A$-field $E$. We
note by deg$_{\tau }$ the degree in indeterminate $\tau $.

\begin{Definition}
An element of $E\{\tau \}$ is called separable, if the
coefficient of its constant element is no null. It called purely
inseparable if it is on the form $\lambda \tau ^{n}$, $n>1$ and
$\lambda \in E$, $\lambda \neq 0$ .
\end{Definition}

Let $H$ be a global field of characteristic $p>0$, and let
$\infty $ one place (a Prime ideal ) of $H$, we will note by
$H_{\infty }$ the completude of $H$ at the place  $\infty $. We
define the degree of function over $A$ by :

\begin{Definition}
Let $a$ $\in A,$ deg $a=\dim _{F_{q}}\frac{A}{aA}$ if $a$ $\neq 0$
and $\deg 0=$ $-\infty $.
\end{Definition}

We extend  deg to $\ K$ by putting deg $x=$ deg $a-$deg $b$ if $0\neq x$ $%
= $ $\frac{a}{b}$ $\in K$.  \\
If $A=\mathbf{F}_{q}[T],$ then the degree function is the usual
polynomial  degree.  Let $Q$ be a no
null ideal  of  $A$, we define the ideal degrees of $%
Q$, noted deg $Q$, by :
\begin{equation*}
\text{deg}Q=\dim _{F_{q}}\frac{A}{Q}.
\end{equation*}

\begin{Lemma}
there exists a rational number  $r$ such that  :
\begin{equation*}
\deg _{\tau }\Phi _{a}=r\deg a\text{.}
\end{equation*}
\end{Lemma}
\begin{proof}
It is easy to see that $\Phi $ is an injection, otherwise since
$K\{\tau \}$ is an integre ring, Ker $\Phi $ is a prime ideal non
null, so maximal in $A$ and by consequence  Im $\Phi $ is a
field, so $\Phi =\gamma $. Since $-$deg$_{\tau }$ define  a no
trivial valuation over Frac$(\Phi (A))$ ( the fractions field of
$\Phi (A)$ ) which is isomorphic to $K$, so $-$deg$_{\tau }$ and
$-$deg are  equivalent valuations over $K$. There is rational
number $r>0 $, such that :
\begin{equation*}
r\ \text{deg}=\text{deg}_{\tau }.
\end{equation*}
\end{proof}
\begin{Corollary}
Let $\Phi :A\mapsto E\{\tau \}$ be an  $A$-Drinfeld module, so
$\Phi $ is injective.
\end{Corollary}

\begin{proposition}
The number  $r$ is a positive integer.
\end{proposition}

\begin{Definition}
The number $r$ is called the rank of the Drinfeld A-module $\Phi.$
\end{Definition}

For example if $A=\mathbf{F}_{q}[T],$ an  A-Drinfeld module of
rank  $r$ is on the  form  : $\ \Phi (T)=a_{1}+a_{2}\tau
+...+a_{r}\tau ^{r}$, where  $\ a_{i}\in E$, $1\leq i\leq r-1$ and
$a_{r}\in E^{\ast }$.

In this case  char $E\mathbf{=}$ $P\neq (0)$ we can define the
notion of height of a Drinfeld module $\Phi$.

For this, we suppose that char$(E)=P\neq 0$. We put  $\
v_{P}:K\mapsto
\mathbb{Z},$ an associate normalized valuation at $P$, this means, if  $a$ $%
\in K$ has a root over  $P$ of order $t$, we have $v_{P}(a)=t$.

For every $a$ $\in A,$ let $w(a)$ the most small integer $t>0$,
where $\tau ^{t}$ occurs  at $\Phi _{a}$ with a non null
coefficient.

\begin{Lemma}
There exists a rational number $h$ such that :%
\begin{equation*}
w(a)=hv_{P}(a)\deg P\text{.}
\end{equation*}
\end{Lemma}

\begin{proposition}
The number $h$ is a positive integer.
\end{proposition}

\begin{Definition}
The integer $h$ is called the height of $\Phi .$
\end{Definition}

For example if $A=\mathbf{F}_{q}[T],$ a Drinfeld $A$-module of height $%
h $ of rank  $r$ is of the form : $\Phi (T)=a_{0}+a_{h}\tau
^{h}+...+a_{r}\tau ^{r},$ where $\ a_{i}\in E$, $0\leq i\leq r-1$
and $a_{r}\in E^{\ast }$.
\subsection{ Endomorphism Ring of a Drinfeld module }
 Let $E$ be an $A$-field and let  $\overline{E}$ a fix algebraic closure.
 Let $\Phi $ and $\Psi $ two Drinfeld  $A$-modules over $E$ of
rank $r>0$ . We define a morphism $\Phi $ to $\Psi $ over $E$ by :

\begin{Definition}
 Let $\Phi $ and $\Psi $ two Drinfeld modules over an $A$-field
\ $L$. A  morphisme of $\Phi $ to $\Psi $ over $E$ is an element $
p(\tau )\in L\mathit{\{\tau \}}$ such that :
\begin{equation*}
p\text{ }\Phi _{a}=\Psi _{a}\text{ }p,\forall a\in A\text{.}
\end{equation*}
\end{Definition}

A no null morphisme is called an isogeny. We note that is
possible only between two Drinfeld modules of same rank.

A invertible isogeny $u$  ( i.e : deg$_{\tau }u=0$) is called a
isomorphism and the module became isomorphic. The set of the
morphims form an $A-$module noted by Hom$_{E}(\Phi ,\Psi ).$

We can see Hom$_{E}(\Phi ,\Psi )$ by the fact that :  a morhism  ( or $E$-morphisme) $p$: $%
\Phi \mapsto \Psi $ from $\Phi $ to $\Psi $ is a morphism of  $A$-modules $%
p$: $(E$,$\Phi $)$\mapsto (L$,$\Psi )$ where $(E$,$\Phi )$ (
respectively  $(E$,$\Psi )$) is $E$ with the structure of
$A$-module given by $\Phi $ (resp $\Psi $).

such morphism is also a morphism of additive groups of  $E$. Then
 $E\{\tau \}$ is a finie  $E[F]-$module so $E\{\tau \}$ is integral over $E[F].$
 So :
\begin{equation*}
E(\tau )=E\{\tau \}\otimes _{E[F]}E(F)=L\{\tau \}\otimes _{A}K,
\end{equation*}

The fraction ring of  $E\{\tau \}$ is noted  $E(\tau )$ ( $E(\tau
)$ is a no commutative field called left field of fractions of
$E\{\tau \}).$

In particularly if $\Phi $ =$\Psi $, the  $E$-endomorphisme
ring   (End$_{E}\Phi =$Hom$_{E}(\Phi ,\Phi $)is a subring of  \
$E\{\tau \}$ and an  $A$-module contained  $\Phi (A)$:

\begin{equation*}
\text{End}_{E}\Phi =\{u\in E\{\tau \}/\forall a\in A,u\Phi
_{a}=\Psi _{a}u\}.
\end{equation*}

Since $\Phi $ is an  injection, $\Phi $ can be naturally
prolonged to an injection \

$\Phi $: $K\mapsto E(\tau ).$ By this  injection we identify on in
$E(\tau ),$ $A$ and $\Phi (A)$ of same $K$ and $\Phi (K).$

Let $F$ the  Frobenius of $E$ we have : $\Phi $ $(A)\subset $ \
End$_{E}\Phi ,$ $ \ F\in $ End$_{E}\Phi .$

\begin{Definition}
Let $\Phi $ and $\Psi $ two Drinfeld $A$-modules over an
$A$-field $E$ and $p$ an isogeny over  $E$ from $\Phi $ to $\Psi
.$

\begin{enumerate}
\item We say that $p$ is separable if and only if  $p(\tau )$ is
separable.
\item We say that  $p$ is purely no separable if and only if $p(\tau
)=\tau ^{j}$ for one $j>0$.
\end{enumerate}
\end{Definition}
\subsection{Norm of Isogeny}

\begin{Definition}
Let $F$ an integer over a ring $A,$ with fractions field $K$. we
note by $N_{K/K(F)}$ the determinant of the $K$-linearly
application of multiplication by  $F$ to $K(F)$ ( it is the usual norm if the extension $%
K(F)/K$ is  separable.
\end{Definition}

We can see that there is a morphism  $N_{K/K(F)}:I_{\overline{A}
}\rightarrow I_{A\text{ }}$ from the fractional ideals groups of  $\overline{%
A}$ to functionary ideals group of $A,$ by this morphism we have:

\begin{proposition}
The norm of an isogeny is a principal ideal.
\end{proposition}

\begin{proposition}
Let $M_{_{fin}}(A)$ the category of primes ideals of $A$ and let
$D(A)$ the  monoïde of ideal of  $A$. There exists a unique
function  :
\begin{equation*}
\chi :M_{_{fin}}(A)\mapsto D(A),
\end{equation*}
multiplicative over the exact sequence and such that $\chi (0)$
$=1$ and $\chi (A/\wp $ $)=\wp $ for every prime ideal $\wp $ of
$A$.
\end{proposition}

\begin{Definition}
The function $\chi $ is called the Euler-Poincare characteristics.
\end{Definition}

We can see $\chi (E^{\Phi })$ and we note it by  $\chi _{\Phi }.$

\begin{proposition}
The ideals $\chi _{\Phi }$ and \ $P^{m}$ are principals  ( in
$A$), and more clearly  $\chi _{\Phi }=(P_{\Phi }(1))$ and
$P^{m}=P_{\Phi }(0)$.
\end{proposition}

\begin{enumerate}
\item We know that the norm of an isogeny is a principal ideal, indeed  $%
N(F)=$ $P_{\Phi }(0)$ and $N(1-F)=$ $(P_{\Phi }(1))$ since $F$ and
$1-F$ are a  $K$-isogenys.

\item We can call  $\chi _{\Phi }$ the divisor of  $E$-points,
this divisor is  analogue at the number of $E$-points for elliptic
curves.

\item $\chi _{\Phi }$ is the annulator of   $A$-module $E^{\Phi }$. We can deduct that  : $E^{\Phi }\subset (\frac{A}{\chi _{\Phi }})^{r}.$

\item The structure of $A$-module $E^{\Phi }$ is stable by the Frobenius endomorphism $F$.
\end{enumerate}
\begin{Corollary}
If there are a Drinfeld  $A$-module  $\Phi ,$ over a field  $E,$
of characteristic $P$ and of degree  $m$ over $A/P$, then the
ideal $P^{m}$ is a principal ideal.
\end{Corollary}

\begin{Remark}
The above Corollary shows that there exists  a  restriction of
the existence of Drinfeld $A$ -modules.
\end{Remark}

\section{Drinfeld Modules Over Finite Fields}

We substitute the Field $E$, by $L$ a finite extension of degree
$n$ of the finite fields $\mathbf{F}_{q}$. The Frobenius  $F$ of
$L$ is \ $F=\tau ^{n}$, so $\mathbf{F}_{q}[F]$ is the center of
$L\{\tau \}$. We put  $m=[L:A/P]$ and $d=$ deg $P$, then $n=m.d$.
The function $-$deg define a  valuation over $K$, the field of
fractions of $A$. Let $\tau :x\mapsto x^{q}$ the  Frobenius of
$\mathbf{F}_{q}$ and let $L$ a finite extension of
$\mathbf{F}_{q}$. We put  $r_{1}=[K(F):K]$ and $r_{2}^{2}$ the
degree of left field  End$_{L}\Phi \otimes _{A}K$ over this
center  $K(F)$.

So a Drinfeld $A$-modules $\Phi $ over $L$ give a  structure of
$A$-module over the additive finite group $L$, this  structure
will be noted  $L^{\Phi }$. Let $\gamma $ the application of $A$
to $L$ which an element $a$ for $A$ associate the constant of
$\Phi _{a}$, then it is easy to see that  $\gamma $ is a ring
homomorphism, and that $\Phi $ and
 $\gamma $ are equal over  $A^{\ast }=\mathbf{F}_{q}^{\ast }$ the set of reversible elements
 of $A$.

\begin{Definition}
Let $\Phi $ be a Drinfeld $A$-module over a finite field $L$. We note by  $%
M_{\Phi }(X)$ the unitary minimal polynomial of $F$ over $K$ .
\end{Definition}

\begin{proposition}
With the above notations  : $M_{\Phi }(X)$ is an element of
$A[X]$, equal to $P_{\Phi }^{\frac{1}{r_{2}}}.$
\end{proposition}

\begin{Corollary}
For two Drinfeld $A$-modules  $\Phi $ and $\Psi $, of rank $r$
over a finite field $L$, then the following are equivalent :

\begin{enumerate}
\item $\Phi $ and $\Psi $ are isogenous,

\item $M_{\Phi }(X)=M_{\Psi }(X),$

\item $P_{\Phi }=P_{\Psi }.$
\end{enumerate}
\end{Corollary}

\begin{proposition}
Let $L$ be a finite extension of degree  $n$ over a finite field
$\mathbf{F}_{q}$, let $F$ the  Frobenius of $L$. Then $L(\tau )$
is a central division algebra over $\mathbf{F}_{q}(F)$ of
dimension $n^{2}.$
\end{proposition}
\begin{Definition}
Every  $u$ $\in L\{\tau \}$ can be writing on this  form $u=\tau
^{h}u^{\prime }$ ( since $L$ is a perfect field) where $u^{\prime
}\in L\{\tau \}$ separable. The integer $h$ is called the height
of  $u$ and will be noted by ht $u.$
\end{Definition}
In the finite fields case, we can see the height of a Drinfeld
$A$-module\ $\Phi $ over finite field $L$, the integer $H_{\Phi }$
as  been :
\begin{equation*}
H_{\Phi }=\frac{1}{r}\inf \{\text{ht }\Phi _{a},0\neq a\in P\}.
\end{equation*}

\begin{Remark}
It is easy to see that  $H_{\Phi }$ is invariant under isogeny
and that
\begin{equation*}
1\leq H_{\Phi }\leq r.
\end{equation*}
\end{Remark}

\begin{proposition}
Let  $\Phi $ be a Drinfeld  $A$-module of rank  $r$ over a finite
field  $L$, the following assertions are equivalent :

\begin{enumerate}
\item There exists a finite extension  $L^{\prime }$ of $L,$ such that
the  endomorphism ring End$_{L^{\prime }}\Phi \otimes _{A}K,$ has
dimension $r^{2}$ over  $K$.

\item Some power of the  Frobenius $F$ of $L$ lies in $A.$

\item $\Phi $ is  supersingular.

\item The field  $K(F)$ has only one prime above $P$.
\end{enumerate}
\end{proposition}
\begin{proposition}
Let $\Phi $ be a Drinfeld  A-module of rank  $r$ and let  $Q$ be
an ideal from $A$ prime with $P,$ then :%
\begin{equation*}
\Phi _{Q}(\overline{L})=(\frac{A}{Q})^{r}\text{.}
\end{equation*}
\end{proposition}
\begin{Corollary}
Then we can deduct that  : $\Phi
_{P}(\overline{L})=(\frac{A}{P})^{r-H_{\Phi }} $.
\end{Corollary}

We can deduct from above mentioned Proposition the following
important result, which characterize the supersingularity :

\begin{proposition}
\ A Drinfeld $A$-module $\Phi $ is supersingular ( $\Phi
_{P}(\overline{L})=0$ ), if and only if, $r=H_{\Phi }$.
\end{proposition}
\begin{Definition}
We say that the field  $L$ is so big if all endomorphism rings
defined over  $\overline{L}$
are also defined over $L$, i.e : End$_{\overline{L%
}}\Phi =$ End$_{L}\Phi .$
\end{Definition}

Two Drinfeld modules $\Phi $ and $\Psi $ are isomorphic, if and
only if, there exists an $a$ $\in $ $L$ such that : $a^{-1}\Phi
_{a}=\Psi _{a}a.$
\begin{Lemma}
Let $\Phi $ be a Drinfeld $ A$-module of rank $r$, over a finite
field $L$, of characteristic $P$. The characteristic polynomial
 of Frobenius endomorphisms $F$ is :%
\begin{equation*}
P_{\Phi }(X)=X^{r}+c_{1}X^{r-1}+...+c_{r-1}X+\mu P^{m},c_{1},...c_{r-1}\in A%
\text{ et }\mu \in \mathbf{F}_{q}^{\ast }\text{.}
\end{equation*}
\end{Lemma}

\begin{Remark}
The fact that constant of the polynomial $P$ is  $\mu P^{m}$ comes
from the fact that $P_{\Phi }(0)=P^{m}$ in $A.$
\end{Remark}

The following Proposition is an analogue of the Riemann's
hypothesis for elliptic curves :

\begin{proposition}
Let $\Phi $ be a Drinfeld $A$-module of rank $r$ over finite
field $L$ which is a finite extension of degree $n$ of \
$\mathbf{F}_{q}$. Then deg $(w)=\frac{n}{r}$\, for every root $w$
of characteristic polynomial $P_{\Phi }(X).$
\end{proposition}

The following result is the Hasse-Weil's analogue for elliptic
curves :
\begin{Corollary}
Let $P_{\Phi }(X)=X^{r}+c_{1}X^{r-1}+...+c_{r}X$ $+\mu P^{m}$
the  characteristic polynomial of a Drinfeld Module $\Phi ,$ of
rank  $r,$ over a finite field $L$. Then:

\begin{equation*}
\forall 1\leq i\leq r-1,\text{deg }c_{i}\leq \frac{i}{r}m\deg P.
\end{equation*}
\end{Corollary}
\begin{proof}
The  proof can be deducted immediately by the above Proposition.
\end{proof}

\section{Drinfeld Modules of rank 2}

In all next of this paper, $\Phi $ will be considered a Drinfeld
$A$-module of rank $2$, And $A=\mathbf{F}_{q}[T]$ for proof and
more details see [1], [12] and [6].

 Our interest for the arithmetic of such modules motive us to
 interesting to their endomorphism rings and their isogenies classes :

\subsection{Endomorphism Ring}

We start by giving the following result which characterized the
order in quadratic extension, to proof see [12] :

\begin{proposition}
Let $O$ be an $A$-order in a quadratic extension $K(F)$, and let
$ O_{K(F)}$ the maximal $A$-order in $K(F)$, then the an $A$-order
$O$ of $ K(F) $ is on the form:
\begin{equation*}
O=A+g.O_{K(F)}\text{,}
\end{equation*}
where $g$ is unitary element of $A$.
\end{proposition}

\begin{Definition}
The element $g$ of $A$ in the  above proposition is called the
conductor of $A$-order $O$.
\end{Definition}

\begin{proposition}
\bigskip Let $\Phi $ be a Drinfeld  $A$-module of rank $2$, over a finite field $L$
of  $m,$ degrees  over $A/P$ and the  Frobenius $F$, and let $P$
be the $ A$-characteristic of  $L$ and $D_{P}$ the completude, at
the place $ P$, of End$_{L}\Phi \otimes _{A}K$.
\begin{enumerate}
\item If $F$ is on the form  $kP^{\frac{m}{2}}(k\in F_{q}^{\ast }),$ then the ring End$_{L}\Phi $
can be identified with a maximal  $A$-order in $_{D_{P}}$,
conversely every maximal  $D_{P}$ can be obtained by this way.

\item Otherwise, the ring End$_{L}\Phi $ can be identified with an
$A$-order in an imaginary quadratic field  $K(F)$. An $A$-order
$O$ of $K(F)$ occur on this way if and only if   \ $F$ $\in $ $O$
, and the conductor of $O$ is prime with  $P$ in two following
cases :

\begin{enumerate}
\item $F$ is on the form$\sqrt{\mu P^{m}}$ with $\mu \in F_{q}^{\ast }$
if $m$ is odd and $\sqrt{\mu P^{m}}$ is imaginary quadratic,

\item $F$ is in the  form $\frac{k_{2}}{k}P^{\frac{m}{2}}$ and $m$ is even and
deg$P$ is odd.
\end{enumerate}
\end{enumerate}
\end{proposition}
\begin{Corollary}
If the conductor of $O$ is prime with $P$, then $O$ is a maximal
$A$ -order in the algebra End$_{L}\Phi \otimes _{A}K$ .
\end{Corollary}

 For the Drinfeld  modules of rank  2, we can specify  End$_{L}\Phi \otimes _{A}K$
 as been in the ordinary case equal to $K(F).$

\begin{proposition}
Let $\Phi $ be a Drinfeld $A$-module of rank $2$ over $L$:

\begin{enumerate}
\item Let $\Phi $ is supersingulary module,

\item Let End$_{L}\Phi \otimes _{A}K=K(F)$.
\end{enumerate}
\end{proposition}

\begin{proof}
Let $r_{1}=[K(F):K]$ and $r_{2}^{2}$ is the degree of left field  End$%
_{L}\Phi \otimes _{A}K$ over this center $K(F)$. Since
$2=r_{1}.r_{1}$, we have two cases ($r_{1}=1$ and $r_{2}=2)$ or
($r_{1}=2$ and $r_{2}=1)$, then in the where case($r_{1}=1$ and
$r_{2}=2$) we have a supersingulary Drinfeld since :
$r_{1}=[K(F):K]=1$ and :
\begin{equation*}
F\in \overline{A^{K(F)}}=\overline{A^{K}}=A.
\end{equation*}%
Otherwise (i.e : $r_{1}=2$ and $r_{2}=1),$ we will have
End$_{L}\Phi \otimes
_{A}K=K(F)$ and End$_{L}\Phi $ is $A$-order in the quadratic field $%
K(F)$.
\end{proof}

\begin{Remark}
 The above result show that the ordinary case and since
  End$_{L}\Phi $ is $A$-order contained
$A[F]$, so to study the maximality of the endomorphism ring
End$_{L}\Phi $,  we can satisfy by study the existence of the
$A$-order containing  \ $A[F]$ and contained in the maximal
$A$-order $O_{K(F)}$ of the algebra $K(F)$.
\end{Remark}

Now we have a no supersingulary Drinfeld A-module and $O_{K(F)}$
is the maximal  $A$-order in the algebra $K(F)$, we interest to
know : when there is an $A$-order $O$ such :

\begin{equation*}
A[F]\subset O\subset O_{K(F)}\text{?}
\end{equation*}

To answer to the above question we have the following result:

\begin{proposition}
Let $\Delta $ $=c^{2}-4\mu P^{m},$ the discriminat of$P_{F}$, the
characteristic polynomial of $F$ the of the finite field $L$,
which is : $ P_{F}(X)=X^{2}-cX+\mu P^{m}$, \ and let  $O_{K(F)}$
the maximally  $A$-order of the algebra $K(F)$.
\begin{enumerate}
\item For every $g\in A$ such that  $\Delta =g^{2}$ $.\omega $, there is a Drinfeld $A$
-module $\Phi $ over $L$ of rank  2, such that :
$O_{K(F})=A[\sqrt{ \omega }]$ and
\begin{equation*}
\text{End}_{L}\Phi =A+g.A[\sqrt{\omega }]\text{.}
\end{equation*}

\item there is not a polynomial  $g$ of $A$ such that $g^{2}$ divide $
\Delta $, there is a ordinary Drinfeld $A$-module $\Phi $ over $L$
of rank 2, such that :
\begin{equation*}
\text{End}_{L}\Phi =O_{K(F)}\text{.}
\end{equation*}
\end{enumerate}
\end{proposition}

\begin{proof}
1) We suppose that there is  $g\in A$, \ such that$\ \Delta $ $=$
$\ g^{2}.\omega $, where $F$ is a root of characteristic
polynomial  $P_{\Phi }$, then we can put : $F=-c/2+\sqrt{\Delta
}/2$ $ =-c/2+g.\sqrt{\omega }/2$, then $A[F]=A[-c/2+g.\sqrt{\omega
}/2]=$ $A[$ $g \sqrt{\omega }/2]\sqsubseteq A+g.A[\sqrt{\omega
}]$ and it is easy to see that in this case the $A$-order
$O_{K(F})$ $=$ $A[\sqrt{\omega }]$ is a maximal $A$-order, and by
the  proposition 4.2, there is a Drinfeld $A$ -module  $\Phi $
such that: End$_{L}\Phi =A+g.O_{K(F)}$.
\\ 2) In the other case :
there is not a polynomial $g$ $\in A$ such that  $g^{2}\mid
\Delta $, and by the  proposition 4.2, there is a Drinfeld
$A$-module $\Phi $ such that the  $A$-order End$_{L}\Phi $ can
not be written on the form  $ A+g.O_{K(F)}$ and in this case it
will be certainly equal to $O_{K(F)}$ .
\end{proof}
\subsection{ Isogenies Classes}

Let $\Phi $ be a Drinfeld $A$-module of rank  2 over a finite
field $L$ with the characteristic polynomial  $P$, and let $m$ =
deg $P$. The characteristic polynomial  \ $P_{\Phi }$ can be
given by the unitary minimal polynomial of  $F$ in $A[X]$,
$M_{\Phi }$, and with the relation  $P_{\Phi }=$ $M_{\Phi
}^{r_{2}}$, $\ r_{2}$ is a root of degree of the left field
End$_{L}\Phi \otimes _{A}K$ over the center $K(F)$.

Let $\overline{K}$ be an algebraic closure of $K$, and $\infty $
a place of $K$ which divide  $\frac{1}{T}$, and let $K_{\infty
}=F_{q}((\frac{ 1}{T}))$, and $ \mathbb{C}_{\infty }$ the
completude of the algebraic closure of $K_{\infty }$.

We fix a plongement $\overset{\_}{K}$ $\hookrightarrow
\mathbb{C}_{\infty }$.

For every $\alpha \in \mathbb{C}_{\infty }$, $\mid \alpha \mid
_{\infty }$ is normalized valuation of  $ \alpha $ ( $\mid
\frac{1}{T}\mid _{\infty }=\frac{1}{q}$).

Let $\theta \in $ $\overset{\_}{K}$, we say that  $\theta $ is an
ordinary number of:

\begin{enumerate}
\item $\theta $ is integral over $A$,

\item $\mid \theta \mid _{\infty }=q^{\frac{md}{2}}$

\item $K$($\theta )/K$ is imaginary, and $[K$($\theta ),K\mathbf{]=}2;$

\item there is only one place of $K$($\theta )$ which divide $\theta $ and Tr$
_{K(\theta )/K}(\theta )\neq 0(P).$
\end{enumerate}

Let $\theta $ an ordinary Weil number $\forall \sigma \in $ col($
\overset{-}{K}\mathbf{/}K$), $\theta ^{\sigma \text{ }}$. We will
noted by  W$^{ord}$ the set of conjugacy classes of an ordinary
Weil numbers of rank 2. Then:

\begin{Theorem}
There is a bijection between  W$^{ord}$ and isogeny classes of
ordinary Drinfeld  $A$-modules of rank $2$ defined over $L$.
\end{Theorem}

To prove the theorem above : let $\theta $ an ordinary Weil
number, we put:

\begin{equation*}
P(x)=Hr(\theta, K\mathbf{;}x\mathbf{)}\text{.}
\end{equation*}

 Then by  (1), (2), (3) and (4) :

\begin{equation*}
P(x)=x^{2}-cx+\mu P^{m}\text{,}
\end{equation*}

 with $\mu \in \mathbf{F}_{q}^{\ast }$ and $c$ $\in A$, $c\neq 0(P)$,
and deg$_{T}$ $c\leq \frac{md}{2}$.

We put $\Gamma =\{c$ $\in A$, $c\neq 0(P)$, deg$_{T}$ $c\leq
\frac{md}{2}\}$.

\begin{Lemma}
Let  $\mu \in \mathbf{F}_{q}^{\ast }$, $c$ $\in \Gamma $, and $E$
the filed of decomposition of $P(x)=x^{2}-cx+\mu P^{m}$ over $K$.
Let $\theta $ a root of $P(x)$. Then $\theta $ verify (1), (2),
and (4) and $[K(\theta ),K\mathbf{]=}2$.
\end{Lemma}

\begin{proof}
Let $B$ be the integral closure of  $A$ in $E$. We suppose that
there is  $ \theta $ a root of  $P(x)$ with $\theta \in B^{\ast
}$. Since the  constant coefficient of $P(x)$ is \ $\mu P^{m}$,
we have  : $\theta \in \mathbf{F} _{q}^{\ast }$. Then: $v_{\infty
}(\theta ^{2}-c\theta )=-$deg$_{T}$ $ c $
> $-md,$ et $\theta ^{2}-c\theta =-\mu P^{m}$ where the contradiction.
We fix then  $\theta $ a root of $P(x)$, we have $\theta \notin
B^{\ast }$ and ($\theta -c)\notin B^{\ast }$. Or $\theta (\theta
-c)=-\mu P^{m}$. Since  $c\neq 0(P)$. There exists exactly  two
primes $ \beta _{1},\beta _{2}$ of $B$ above $P$ and $\beta
_{1}\mid \theta ,$ $ \beta _{2}\mid \theta -c$. In particularly
$[E:K]=2$. We work in  $ \mathbb{C}_{\infty },$ we have :
\begin{equation*}
v_{\infty }(\theta )+v_{\infty }(\theta -c)=-md.
\end{equation*}
Since $v_{\infty }(c)=-$deg$_{T}$ $c$ $\geq \frac{-md}{2}$, we
have $v_{\infty }(\theta )<0$. We suppose that $v_{\infty
}(\theta )$ or $v_{\infty }(\theta -c)\neq \frac{-md}{2}$. Even
we replace  $\theta $ by $\theta -c$, we can suppose :
\begin{equation*}
v_{\infty }(\theta )<\frac{-md}{2}\text{.}
\end{equation*}
Then :
\begin{equation*}
v_{\infty }(\theta -c)=\inf (v_{\infty }(\theta ),v_{\infty
}(c))=v_{\infty }(\theta )\text{.}
\end{equation*}%
Where from the contradiction.
\end{proof}
\begin{Corollary}

1) Let $\mu \in \mathbf{F}_{q}^{\ast }$ and $c$ $\in \Gamma $,
then if  $\theta $ is a root of $x^{2}-cx+\mu P^{m}$, $\theta $
is a ordinary Weil number, if and only if $K(\theta )/K$ is
imaginary.

2)If $md\equiv 1(2)$, then $ \forall \mu \in \mathbf{F}_{q}^{\ast
}$ and $\forall c$ $\in \Gamma $, the roots of $x^{2}-cx+\mu
P^{m}$ are the Weil numbers.
\end{Corollary}

To simplify our speech, we suppose  $p\neq 2$. And we put
$md\equiv 0(2) $.

\begin{Lemma}
Let $\mu \in \mathbf{F}_{q}^{\ast }$ and $c$ $\in \Gamma $ with deg$_{T}$ $%
c\leq \frac{md}{2}$. Let $\theta $ a root of $x^{2}-cx+\mu P^{m}$.
Then $\theta $ is a Weil number, if and only if, $-\mu \notin (%
\mathbf{F}_{q}^{\ast })^{2}$.
\end{Lemma}

\begin{proof}
By the  Hensel lemma : $P^{m}\in (K_{\infty }^{\ast })^{2}$. We
have :

\begin{equation*}
v_{\infty }(\frac{c}{\sqrt{P^{m}}}\ ) =\frac{md}{2}\ -deg_{T}\ c>0
\text{.}
\end{equation*}
Or $\frac{\theta }{\sqrt{P^{m}}}\ \ $ is a root:
\begin{equation*}
\ x^{2}-\frac{c}{\sqrt{P^{m}}}x+\mu =0\text{.}
\end{equation*}%
Or : $x^{2}-\frac{c}{\sqrt{P^{m}}}x+\mu \equiv x^{2}+\mu $ \ \
($\frac{1}{T} \mathbf{F}_{q}^{{}}[[\frac{1}{T}]])$.

By the Hensel lemma $\theta \notin (\mathbf{F}_{q}^{\ast
})^{2}\Leftrightarrow -\mu \notin (\mathbf{F}_{q}^{\ast })^{2}$.
\end{proof}

\begin{Lemma}
Let $\mu \in \mathbf{F}_{q}^{\ast }$ and $c$ $\in \Gamma $ with
deg$_{T}$ $c= \frac{md}{2}$, and we note  $c_{0}$ the term of
more high degree of  $c$. We suppose  $c_{0}^{2}\neq -4\mu $. Let
 $\theta $ a root of $x^{2}-cx+\mu
P^{m}$. Then $\theta $ is a Weil number, if and only if, $%
x^{2}-c_{0}x+\mu $ is irreducible in  $\mathbf{F}_{q}[X]$.
\end{Lemma}
\begin{proof}
In this time we choose,  $\sqrt{P^{m}}$ such that :
\begin{equation*}
\sqrt{P^{m}}(\frac{1}{T})^{\frac{md}{2}}\equiv 1(\frac{1}{T})
\end{equation*}%
Then : $\frac{c}{\sqrt{P^{m}}}\equiv 0(\frac{1}{T})$. Then :
\begin{equation*}
x^{2}-\frac{c}{\sqrt{P^{m}}}x+\mu \equiv x^{2}-c_{0}x+\mu \ \
(\frac{1}{T}).
\end{equation*}%
We apply the Hensel lemma since $x^{2}-c_{0}x+\mu $ have two
simple  roots.
\end{proof}

If $c_{0}^{2}=-4\mu $, we put $\Delta =c^{2}-4\mu $. Then $\theta
$ is a Weil number if and only if deg $_{T} \Delta \cong 1$(2). If
\ deg $_{T} \Delta \cong 0$ (2)and the term of more high degree in
$\Delta $ is not a square in $\mathbf{F}_{q}^{\ast}$.

We note some remarks about the Weil numbers:

\begin{enumerate}
\item The Weil numbers are defined for all the rank,


\item The  Frobenius $F$ over $L$ is a Weil number.
\end{enumerate}

We know, By [7] and [12] that the characteristic polynomial
$P_{\Phi }$ \ of a Drinfeld $A$-module of rank  $2$, is one of
the four following form  :

\begin{proposition}
Let $\Phi $ be a Drinfeld  $A$-module of rank $2$ over the finite
field $L= \mathbf{F}_{q^{n}}$ and let  $P$ the characteristic of
$L$. We put  $ m=[L:A/P]$ and $d=$deg $P$. The characteristic
polynomial  $P_{\Phi }$ is on the form  :

1)$P_{\Phi }(X)=X^{2}-cX+\mu P^{m}$, where $c^{2}-4\mu P^{m} $ is
imaginary, $c\in A$, $(c,P)=1$ and $\mu \in \mathbf{F}_{q}^{\ast
}$, if $\Phi $ is an ordinary module.\\
 And the characteristic polynomial  $P_{\Phi }$ for a supersinglary case  is on the form :

2)$P_{\Phi }(X)=X^{2}+\mu P^{m}$, with $\mu \in
\mathbf{F}_{q}^{\ast }$ , if $m$ is odd,

3)$\ \ P_{\Phi }(X)=X^{2}+c_{0}X+\mu P^{m},$ if $m$ is even and
$d=\deg P$ is odd, $\mu \in \mathbf{F}_{q}^{\ast }$ and $c_{0}\in
\mathbf{F}_{q}$.

4)$P_{\Phi }(X)=(X+\mu P^{\frac{m}{2}})^{2}$ if $m$ is even.
\end{proposition}


We can resumed  all the cases in following manner :

\begin{enumerate}
\item For the ordinary case, the characteristic polynomial is on
the form :%
\begin{equation*}
P_{\Phi }(X)=X^{2}-cX+\mu P^{m}\text{,}
\end{equation*}%
\ such that : $2$ deg $c<\deg P.m$ or $2\deg c=$deg $P.m$ and
$X^{2}-a_{0}X+\mu $ is irreducible over $\mathbf{F}_{q^{n}}$
where  $a_{0}$ is the coefficient of more big degree of $c$. For
the supersingular case, we have the two following cases :

\item $\deg P$ is even or $-\mu \notin (\mathbf{F}_{q}^{\ast })^{2}$.

\item $X^{2}+c_{0}X+\mu $ is irreducible over $\mathbf{F}_{q}$.
\end{enumerate}

Then we can now calculate the number of characteristic polynomial
which corresponding of the number of isogeny classes :

\begin{Lemma}
\begin{equation*}
\#\{\text{classes d'isog\'{e}nies}\}=\#\{P_{\Phi }\}.
\end{equation*}
\end{Lemma}

We start by calculate the ordinary case. The case(1) give us the
generally number  $c,$ which  corresponding the number of isogeny
classes for an ordinary Drinfeld module, and by using the 2,3 and
4 we find the number of isogeny classes for an supersingulary
Drinfeld module. Indeed in this case  1, the principal condition
la condition that we have above $c$, other than this prime with
$P,$ is the  Riemann condition which certify that $c^{2}-4\mu
P^{m}$ is imaginary which condition can be wrote by : $\deg c\leq
\frac{m.d}{2}.$

We distinguish, two cases :

1) The case where the number  \ \ \ $m.d$ is odd, that means  $\
m $ and \ $d$ are odds. We will have $q^{[\frac{m}{2}d]+1}$
polynomials of degree less or equal than $[\frac{m}{2}.d]$ where  $%
[..]$ is the  partial entire  ). In next, we eliminate the
polynomials $c$ which are not primes with  $P$, that means
divisible by $P$. We can remark that for each  $c$ divisible by
$P$, there is a polynomial $Q$ such that  $c=Q.P$ then the
cardinal of such polynomials  $c$ which are  divisible by $P$ is
equal to the cardinal of the set of $Q$ which is the order of
$q^{\frac{m-2}{2}d+1}( since $ deg $Q$ $\leq \frac{m-2}{2}d)$. If
we consider than  $\mu \in \mathbf{F}_{q}^{\ast }$, we will have :

$\#\{P_{\Phi },\Phi :$ordinary(1)$\}=(q-1)(q^{[\frac{m}{2}d]+1}-q^{[\frac{%
m-2}{2}d]+1})$.

2) For the case where the number $\frac{m}{2}.d$ \ is even that
means that at least one of the $m$ or $d$ is even, we will
exclude the minimal polynomial associated to the corresponding
modules which are not irreducibles and the condition over $c$
became then :

deg$c < m.\frac{d}{2}$ and the polynomial $X^{2}-c_{0}X+c$ is
irreducible where  $c_{0}$ is the  coefficient of more big degree
of $c$, with the prime condition of $c$ and $P$ we will have in
this case:

 $\#\{P_{\Phi }$, $\Phi $ \ ordinary ( 1)$\}=(q-1)(\frac{(q-1)}{2}%
q^{\frac{m}{2}d}$ $-q^{\frac{m-2}{2}d+1}).$

For the case where the characteristic polynomial is on the form

$P_{\Phi }(X)=X^{2}+\mu P^{m}$, where $\mu \in F_{q^{n}}$ if $m$
is even. We will have  $q-1$ possibilities, and $q^{2}-q$
possibilities for the case 3, and finally we will have  $q-1$
possibilities for the case  4. So we can calculate the  cardinal
of the isogeny  classes  of a Drinfeld  module of rank  $2$:
\begin{proposition}
 Let $\Phi $ a Drinfeld  $A$-module of  rank $2$ over a finite field $L=F_{q^{n}}$ and
 let $P$ be the  $A$-characteristic of $L$. We put $m=[L:A/P]$ and $d=$deg $P$ :

\begin{enumerate}
\item  $m$ is odd and $d$ is odd :%
\begin{equation*}
\#\{P_{\Phi },\Phi :\text{ordinary}(1)\}=(q-1)(q^{[\frac{m}{2}d]+1}-q^{[%
\frac{m-2}{2}d]+1}+1).
\end{equation*}

\item $ m$ is even and $d$ is odd:%
\begin{equation*}
\#\{P_{\Phi
}\}=(q-1)[\frac{q-1}{2}q^{\frac{m}{2}d}-q^{\frac{m-2}{2}d+1}+q].
\end{equation*}

\item $m$ is even and $d$ is even:%
\begin{equation*}
\#\{P_{\Phi
}\}=(q-1)[\frac{q-1}{2}q^{\frac{m.}{2}d}-q^{\frac{m-2}{2}d}+1].
\end{equation*}
\end{enumerate}
\end{proposition}

\subsubsection{Characteristic of  Euler-Poincare}

\bigskip Let $\Phi $ be a Drinfeld $A$-module of  rank $2,$ over a finite field
$L=\mathbf{F}_{q^{n}}$ and the polynomial characteristic $P_{\Phi
}.$ We have seen above that  $\chi _{\Phi }=(P_{\Phi }(1))$, this
give us the possibility to deduct that if  $\Psi $ is an other
Drinfeld  $A$-module of rank $2$ over finite field $L$ \ of the
polynomial  characteristic  $ P_{\Psi }$ and the  Euler-Poincare
characteristic $\chi _{\Psi }$, then:

$\ $
\begin{equation*}
\ \chi _{\Phi }=\chi _{\Psi }\Longleftrightarrow \exists \lambda
\in F_{q^{n}}^{\ast }:P_{\Phi }(1)=\lambda P_{\Psi }(1).
\end{equation*}

that means that the cardinal of the set of characteristic of the
Euler-Poincare, can be deduct by the cardinal  of the set of
characteristic polynomial, and we have :
\begin{equation*}
\#\{\chi _{\Phi }\}\leq \frac{\#\{P_{\Phi }\}}{q^{n}-1}\text{.}
\end{equation*}

About the characteristic of the Euler-Poincare we can enumerate
the following remark :
\begin{enumerate}
\item The characteristic of
the Euler-Poincare is the analogue of the number of points of the
elliptic curve over this finite field.
\item For two elliptic curves, it is sufficient to have the same
number of points to be isogenous, but for two Drinfeld modules it
is no sufficient for two Drinfeld modules $\Phi $ and $\Psi $ are
isogenous, since this  two modules are isogenous if and only if
$P_{\Phi }=P_{\Psi }$, or the fact that $\chi _{\Phi }=\chi
_{\Psi }$ implies only that there is $\lambda \in F_{q^{n}}^{\ast
}$ such that :
\begin{equation*}
P_{\Phi }(1)=\lambda P_{\Psi }(1)\text{.}
\end{equation*}
\end{enumerate}

Then we can have a expression for the cardinal of the set of
characteristic of  Euler-Poincare:

\begin{proposition}
 Let $\Phi $ be a Drinfeld A-module of rank  $2$ over the finite field
  $L=\mathbf{F}_{q^{n}}$ and let $P$ be the characteristic of
$L$. We put $m=[L:A/P]$ and $d=$deg $P$. There exists $H,B\in L,$
such that:
\begin{equation*}
\#\{\chi _{\Phi }\}=H+B
\end{equation*}
where $H$ and $B$ verifies:
\begin{equation*}
\#\{P_{\Phi }\}=(q-1)H+(q-2)B.
\end{equation*}
\end{proposition}

\begin{proof}
Let $\Phi $  and $\Psi ,$ two Drinfeld $A$-modules over
$\mathbf{F} _{q^{n}}$, $P_{\Phi }(1)=1-c+\mu P^{m}$ and $P_{\Psi
}(1)=1-c^{\prime }+\mu ^{\prime }P^{m}$. Then $\chi _{\Phi }$
$=\chi _{\Psi },$ if and only if, $\lambda \in F_{q^{n}},$ such
that  $P_{\Phi }(1)=\lambda P_{\Psi }(1) $ then : $1-c+\mu P^{m}$
$=$ $\lambda -\lambda c^{\prime }+\lambda \mu ^{\prime }P^{m}$.
that means that : $\mu =\mu ^{\prime }$ and $\ c^{\prime }=\lambda
^{-1}(1-c+\lambda )$, then the $\lambda $ are of order of $q-2$ (
car $\lambda \in F_{q}-\{0,1\}$) . At past the non prime
condition with  $P$, we will have, if such divisor $ Q$ exists, \
$Q.P=1+\lambda +\lambda c^{\prime }$ and then  deg $Q=-d+\deg
c^{\prime }\leq $ $\frac{(m-2}{2})d$. Then the cardinal of these
$Q$ is equal to the cardinal of these $c^{\prime }$, which is
$q^{[\frac{m}{2}d]+1}-q^{[ \frac{m-2}{2}d]+1}$ which is $B$, and
the couple  $(\lambda ,t^{\prime })$ are of the order
$(q-2)(q^{[\frac{m}{2} d]+1}-q^{\frac{[m-2}{2}d]+1})$ . We can
have  $H$ by the equation  : $\#\{P_{\Phi
}\}=(q-1)H+(q-2)B\Longrightarrow H=\frac{1}{ q-1}(\#\{P_{\Phi
}\}-(q-2)B):$\bigskip \bigskip\ $\ \ \ $ we start by the case :$\
\ \ \ $1) $m$ is odd and $d$ is odd :
\begin{equation*}
H=\frac{1}{q-1}q^{[\frac{m}{2}d]+1}-\frac{1}{q-1}q^{\frac{[m-2}{2}d]+1}+1.
\end{equation*}%
$\bigskip $2) $m$ is even and $d$ is odd:%
\begin{equation*}
H=\frac{1+2q-q^{2}}{2q-2}q^{\frac{m}{2}d}-\frac{1}{q-1}q^{\frac{m-2}{2}d+1}+q%
\text{.}
\end{equation*}%
\ 3) $m$ is even and $d$ is even:%
\begin{equation*}
H=\frac{1+2q-q^{2}}{2q-2}q^{\frac{m}{2}d}-\frac{1}{q-1}q^{\frac{m-2}{2}d+1}+1%
\text{.}
\end{equation*}
\end{proof}

Finally we recuperate the values of $\#\{\chi _{\Phi }\}:$

\begin{proposition}
 Let $\Phi $ be a Drinfeld A-module of rank  $2$ over a finite field $%
L=\mathbf{F}_{q^{n}}$ and let $P$ the  $A$- characteristic of $L$.
We put $m=[L:A/P]$ and $d=$deg $P$ :

\begin{enumerate}
\item $m$ is odd and $d$ is odd :

\begin{equation*}
\#\{\chi _{\Phi }\}=\frac{q}{q-1}q^{[\frac{m}{2}d]+1}-\frac{q}{q-1}q^{[\frac{%
m-2}{2}d]+1}+1
\end{equation*}

\item $m$ is even and $d$ odd:%
\begin{equation*}
\#\{\chi _{\Phi }\}=\frac{q^{2}+1}{2q-2}q^{\frac{m}{2}d}-\frac{q}{q-1}q^{%
\frac{m-2}{2}d+1}+q
\end{equation*}

\item $m$ is even and $d$ is even:%
\begin{equation*}
\#\{\chi _{\Phi }\}=\frac{q^{2}+1}{2q-2}q^{\frac{m}{2}d}-\frac{q}{q-1}q^{%
\frac{m-2}{2}d+1}+1\text{.}
\end{equation*}
\end{enumerate}
\end{proposition}


\begin{thebibliography}{99}
\bibitem{Angeles} Bruno Angles. Th\`{e}se de Doctorat : Modules de Drinfeld
sur les corps finis, Universit\'{e} Paul Sabatier-Toulous III, no
d'ordre 1872, (1994).

\bibitem{B.Angles} Bruno Angles. One Some Subring of Ore Polynomilas
Connected with Finite Drinfeld Modules, J. Algebra 181\textbf{\ }(1996)%
\textbf{, }no.2, 507--522.

\bibitem{Gosse} David Goss. Basic Structures of Function Field Arithmetic,
Volume 35 Ergbnise der Mathematik und ihrer Grenzgebiete,
Springer.

\bibitem{Drinfeld2} V.G. Drinfeld. Modules Elliptiques. Math, USSR Sbornik,
94 (136), 594-627, 656, (1974).

\bibitem{Drinfeld1} V.G. Drinfeld. Modules Elliptiques II Math, USSR
Sbornik, 102 (144), No 2, 182-194,325, (1977).

\bibitem{Gekeler 1} Ernst-Ulrich Gekeler. On Finite Drinfeld Module. J.
Algebra 141, (1991), 187-203.

\bibitem{Gekeler} Ernst-Ulrich Gekeler and Brian A. Snyder. Drinfled Modules
Over Finite Fields . Drinfeld Modules, Modular Schemes and
Application. Alden-biesen, (1996).

\bibitem{Potemine} Igor Potemine. Th\`{e}se de Doctorat : Arithm\'{e}tiques
des Corps Globaux de Fonctions et G\'{e}om\'{e}trie des
Sch\'{e}mas Modulaires de Drinfeld, de l'Universit\'{e} Joseph
Fourier( Grenoble I), (1997)

\bibitem{Silvermann} Joseph. H. Silverman The Arithmetic of Elliptic Curves.
Graduate Texts in Mathematics, 106.

\bibitem{YU1} J.K.Yu. A Classe Number Relation over Functions Fields, J.
Number Theory, 54, (1995), 318--340.

\bibitem{YU} J-K.\ Yu. Isogenies of Drinfeld Modules over Finite Fields, J.
Number Theory 54 (1995), \ no 1, 161--171.

\bibitem{Deuring} M. Deuring. Die \ Typen der Multiplikatorenringe
Ellipticher Funktionenkorper, Abh. Math.sem.Univ.Hamburg, 14 (
1941), 197-272.

\bibitem{Tsfa-Vladut} M. A.Tsfasman-S. G. Vladut. Algebraic-Geometric Codes,
Mathematics and Applications, Dordrecht et al, (1991).

\bibitem{Schoof} R. Shoof. Nonsingular Plane Cubic Curves over Finite
Filelds, Journal of combinatory theory, series A 46, (1987),
183-211.

\bibitem{Reiner} I. Reiner. Maximal Orders. Academic Presse, (1975).

\bibitem{Rok} H.G. Ruck. A Note on Elliptic Curves Over Finite Fields. Math.
Comp. 49, no179, (1987), 301--304.

\bibitem{Waterhouse} W. C. Waterhouse. Abelian Varieties \ Over Finite
Fields. Ann. Sci. Ecole Norm. Sup2, (1969), 521-560.
\end{thebibliography}
\end{document}